\documentclass[11pt]{article}
\usepackage{amsfonts,amsbsy,amssymb,amsmath,graphicx,amsthm,float,subfigure,natbib,color,url,inputenc,geometry}
\usepackage[english]{babel} 
\usepackage{epstopdf}
\usepackage{ifpdf}
\usepackage{verbatim}

\begin{document}

\bibliographystyle{natbib}

\title{\textbf{Response to: ``Limitations of the Method of Lagrangian Descriptors'' [http://arxiv.org/abs/1510.04838]}}
\author{F. Balibrea-Iniesta$^{1}$, J. Curbelo$^{2}$, V. J. Garc\'ia-Garrido$^{1}$, \\ C. Lopesino$^{1}$, A. M. Mancho$^{1}$, C. Mendoza$^{1}$, S. Wiggins$^{3}$\\ \\
$^1$Instituto de Ciencias Matem\'aticas, CSIC-UAM-UC3M-UCM,
\\ C/ Nicol\'as Cabrera 15, Campus Cantoblanco UAM, 28049
Madrid, Spain.\\[.2cm]
$^2$Laboratoire de G\'eologie de Lyon (CNRS/ENS-Lyon/Lyon1),
 Lyon, France.\\[.2cm]
$^3$School of Mathematics, University of Bristol, \\Bristol BS8 1TW, United Kingdom.}

\maketitle

\begin{abstract}
This Response is concerned with the recent  Comment of Ruiz-Herrera, ``Limitations of the Method of Lagrangian Descriptors'' [http://arxiv.org/abs/1510.04838], criticising the method of  Lagrangian Descriptors. In spite of the significant body of literature asserting the contrary, Ruiz-Herrera claims that the method  fails to reveal the presence of  stable and unstable manifolds of hyperbolic trajectories in  incompressible systems and in almost all linear systems. He supports this claim by considering the method of Lagrangian descriptors applied to three specific examples. However in this response we show that Ruiz-Herrera does not understand the proper application and interpretation of the method and, when correctly applied, the method beautifully and unambiguously detects the stable and unstable manifolds of the hyperbolic trajectories in his examples.
\end{abstract}

\section{Introduction}

Before analyzing in detail the examples and  claims of Ruiz-Herrera, we first describe the set-up for Lagrangian Descriptors as originally described in \cite{chaos}, and further developed in \cite{cnsns}.
We consider a vector field,

\begin{equation}
  \frac{d\mathbf{x}}{dt} = \mathbf{v}(\mathbf{x},t), \quad \mathbf{x} \in \mathbb{R}^n \;,\; t \in \mathbb{R}
\end{equation}

\noindent
where $\mathbf{v}(\mathbf{x},t) \in C^r$ ($r \geq 1$) in $\mathbf{x}$ and continuous in time. We denote a trajectory of this vector field by $\mathbf{x}(t;\mathbf{x}_{0})$, with initial condition $\mathbf{x}(t_0;\mathbf{x}_{0}) = \mathbf{x}_{0}$.

The classical definition of the Lagrangian Descriptor (LD, the function $M$), as introduced in \cite{chaos,prl}, is:

\begin{equation}
M(\mathbf{x}_{0},t_{0},\tau) = \displaystyle{\int^{t_{0}+\tau}_{t_{0}-\tau} \|\dot{\mathbf{x}}(t;\mathbf{x}_{0})\| \; dt}
\label{M_classic}
\end{equation}

\noindent
where $ \|\dot{\mathbf{x}}(t;\mathbf{x}_{0})\| \equiv \sqrt{<\dot{\mathbf{x}}(t;\mathbf{x}_{0}) , \dot{\mathbf{x}}(t;\mathbf{x}_{0})>}$. So the LD  
is a function 
depending on the initial condition of a trajectory $\mathbf{x}(t_0;\mathbf{x}_{0}) = \mathbf{x}_{0}$ and  on a time interval $[t_0-\tau, t_0+\tau]$. Hence  \eqref{M_classic} is the arclength of trajectories over the time interval $[t_0-\tau, t_0+\tau]$. Phase space structure, e.g. stable and unstable manifolds of hyperbolic trajectories,  is encoded in the properties of the function $M(\mathbf{x}_{0},t_{0},\tau)$. Concerning this point, in his Comment  Ruiz-Herrera states the following,

\medskip
\noindent
\begin{quotation}
\em The method of Lagrangian descriptors says that the invariant manifolds of saddle points in (2.2) are given by ``singular points'' (i.e. ``non-smooth points'') of the contour lines of $M$, see [10, 4, 7].\footnote{Note that in Ruiz-Herrera Comment, [10] corresponds to \cite{carlos}, [4] corresponds to \cite{cnsns} and [7] corresponds to \cite{prl}.} 
\end{quotation}
\medskip

\noindent
However, {\em no such assertion is made in any the references quoted by Ruiz-Herrera}. In fact, {\em it is not true}, and it explains some of the confusion of Ruiz-Herrera and why, in his examples, he does not  detect the stable and unstable manifolds of hyperbolic trajectories with his incorrect usage of LDs. The term ``singularities'' does {\em not} refer to properties of the contour lines of $M(\mathbf{x}_{0},t_{0},\tau)$,  but to points at which certain directional derivatives  of $M(\mathbf{x}_{0},t_{0},\tau)$ do not exist.  While this is discussed, somewhat, in the references \cite{prl, cnsns}, it is made precise in \cite{carlos}.

Henceforth we will consider  his examples separately (and refer to specific equation numbers in his Comment when necessary).

\section{First Example of Ruiz-Herrera}

In his equation (2.3) Ruiz-Herrera considers the vector field:

\begin{eqnarray*}
\dot{x} & = & \lambda x, \\
\dot{y} & = & -\lambda y, \quad \lambda >0.
\end{eqnarray*}

\noindent
Note that the origin is a hyperbolic fixed point, with the $x$-axis corresponding to its unstable manifold and its $y$-axis corresponding to its stable manifold. In his figure 1 Ruiz-Herrera plots the contours of $M$ (for $\lambda =1$ and $\tau=20$). For this case the contours of  $M$ do appear to converge to the stable and unstable manifolds of the hyperbolic fixed point. He claims (without proof, and incorrectly as we will show) that this is a result of the equal expansion and contraction rates of the saddle.

\begin{figure}[H]
\centering
{\includegraphics[scale = 0.4]{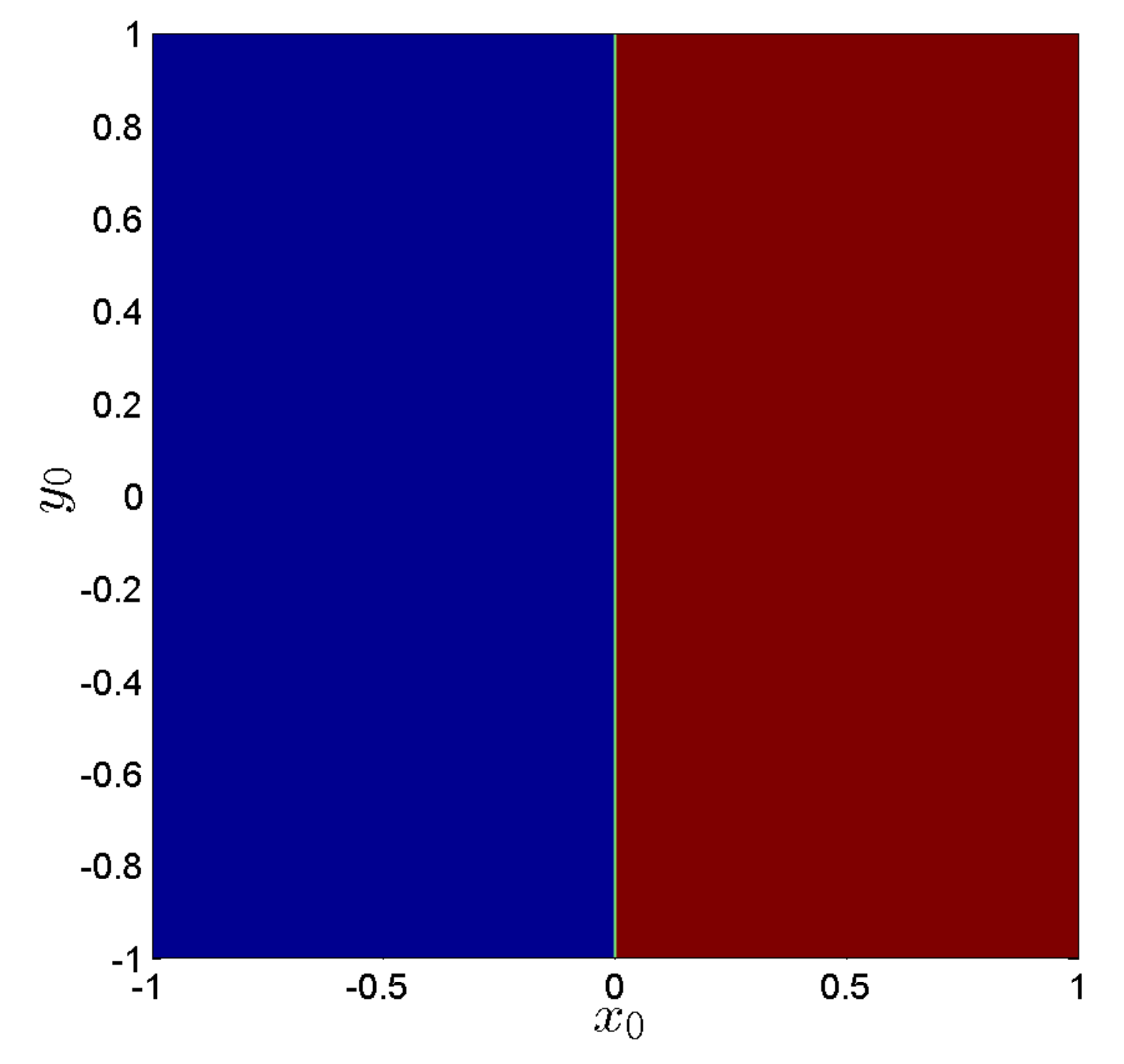}}
{\includegraphics[scale = 0.4]{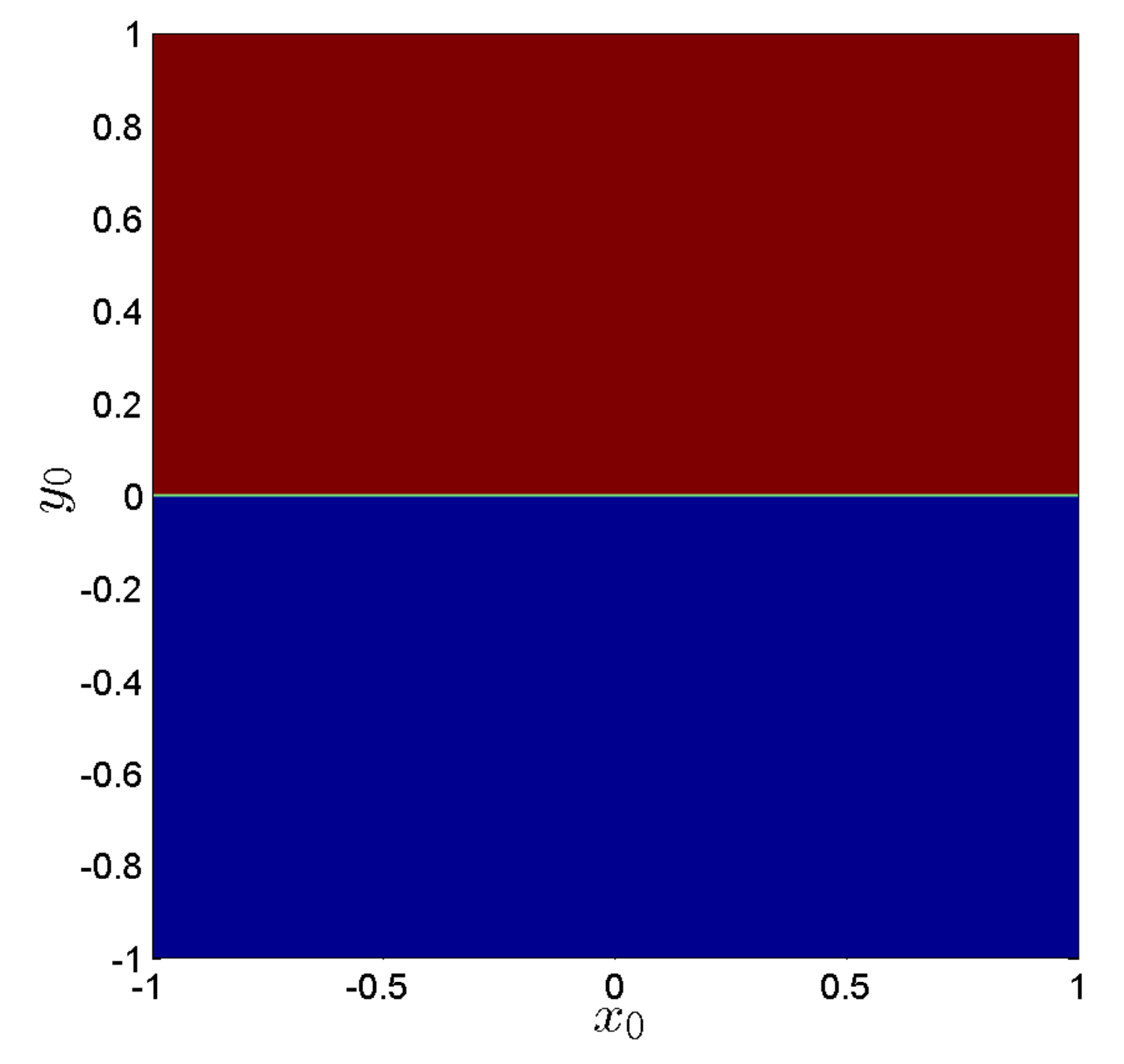}}
\caption{The left-hand panel shows contours of $\partial_{x_{0}} M$ and the right-hand panel shows contours of $\partial_{y_{0}} M$.  This figure should be compared with figure 1 of the Comment of Ruiz-Herrera. }
\label{Afig1}
\end{figure}

From now until the end of the paper, the value of variables $t_{0}$ and $\tau$ will be pre-fixed for every analyzed example. Without loss of generality we can take $t_{0}=0$ since every example considered is an autonomous (time independent) system. As a consequence of this, the domain of the function $M(\textbf{x}_{0},t_{0},\tau )=M(x_{0},y_{0},t_{0},\tau )$ will be the plane $\mathbb{R}^{2}$ of initial conditions $\textbf{x}_{0}=(x_{0},y_{0})$. Now we also compute $M$ for $\lambda =1$ and $\tau =20$ as Ruiz-Herrera in his Comment.
In the left-hand panel of Fig. \ref{Afig1} we show the contours of $\partial_{x_{0}} M = \partial M/\partial x_{0}$, the partial derivative of $M$ with respect to $x_{0}$. Here we can see that $\partial_{x_{0}} M$ changes in sign on the stable manifold $\lbrace x_{0}=0 \rbrace$, exactly as expected.
In the right-hand panel of Fig. \ref{Afig1} we show the contours of $\partial_{y_{0}} M = \partial M/\partial y_{0}$, the partial derivative of $M$ with respect to $y_{0}$.  Similarly, $\partial_{y_{0}} M$ changes in sign on the unstable manifold $\lbrace y_{0}=0 \rbrace$, as expected. These two special features highlight the two manifolds over a region of the hyperbolic point.

\section{Second Example of Ruiz-Herrera}

In his equation (3.5) Ruiz-Herrera considers the vector field:

\begin{eqnarray*}
\dot{x} & = & f(x), \\
\dot{y} & = & - y \hspace{0.1cm} \cdot \hspace{0.1cm} f'(x),
\end{eqnarray*}

\noindent
where $f(x)$ is given in a complicated expression in (3.4), which we will not reproduce here. One feature of the definition of $f(x)$ is that the origin is a hyperbolic fixed point whose stable manifold is given by the $y$-axis and whose unstable manifold is given by the $x$-axis. The other feature is that it allows him to show that the contour lines of $M$ are horizontal lines in a neighborhood of the stable manifold. He claims that this  demonstrates that the method of LDs does not detect the stable and unstable manifolds. However, as in the previous example, we will show  that this is not true by using $M$ properly.

\begin{figure}[H]
\centering
{\includegraphics[scale = 0.4]{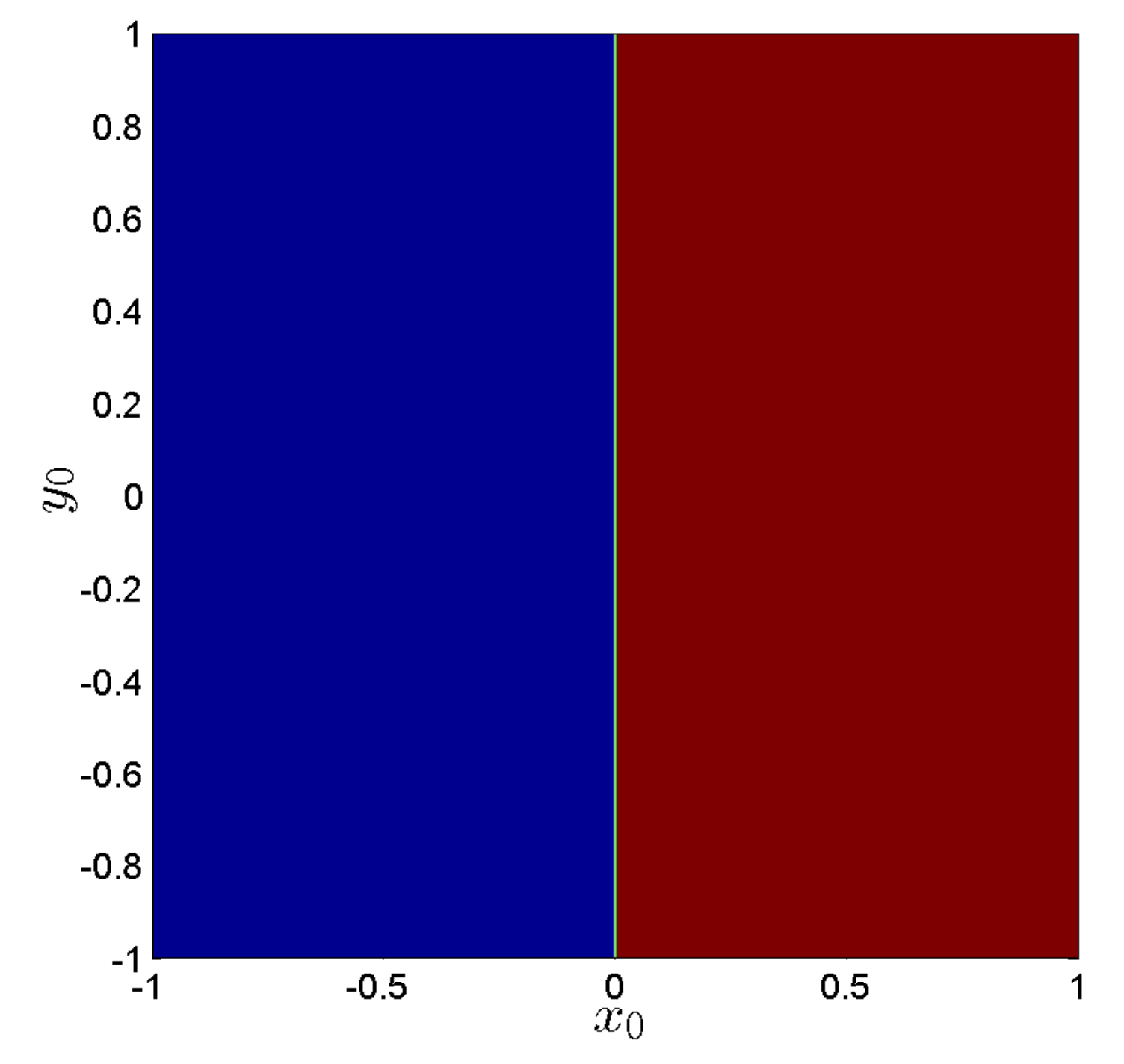}}
{\includegraphics[scale = 0.4]{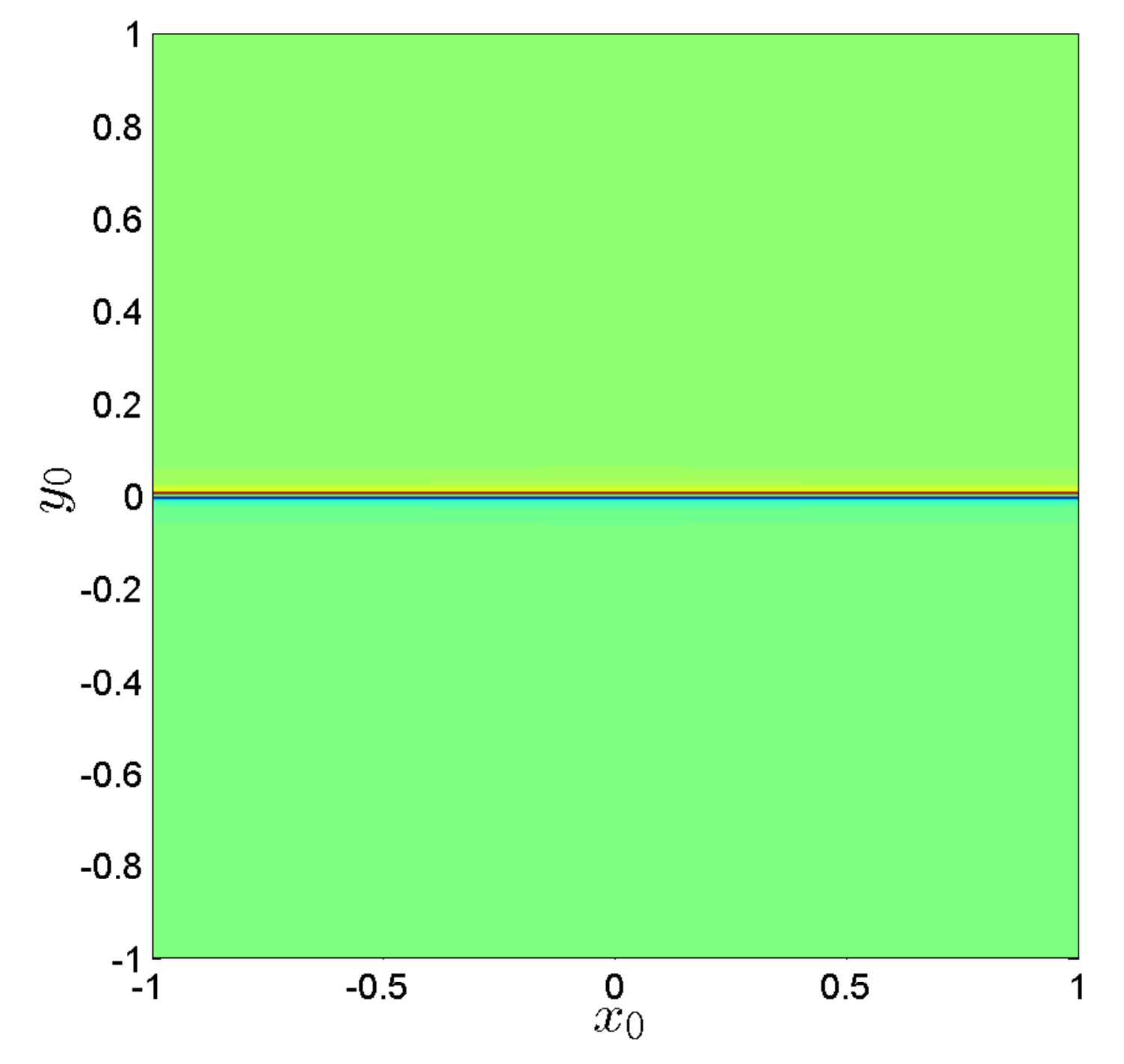}}
\caption{The left-hand panel shows contours of  $\partial_{x_{0}} M$ and the right-hand panel shows contours of $\partial_{y_{0}} M$. This figure should be compared with figure 2 of the Comment of Ruiz-Herrera. }
\label{Afig2}
\end{figure}

We  have computed $M$ for this example using $\tau =10$.
In  the left-hand panel of Fig. \ref{Afig2} we show contours of $\partial_{x_{0}} M$. Its sign changes on the stable manifold, exactly as expected.
In the right-hand panel of Fig. \ref{Afig2} we show contours of $\partial_{y_{0}} M$. Similarly to the left-hand panel, $\partial_{y_{0}} M$ experiences a change in sign on the unstable manifold, although this change is more abrupt than in the previous example and this might be seen only as an anomaly of $\partial_{y_{0}} M$.

We will not consider the third example of Ruiz-Herrera since it is a linear saddle point and illustrates exactly the same type of misunderstanding of LDs as his other examples. We will go directly to his fourth and final  example.

\section{Fourth Example of Ruiz-Herrera}

In his equation (4.9) Ruiz-Herrera considers the vector field,

\begin{eqnarray*}
\dot{x} & = & \lambda x, \\
\dot{y} & = & -\mu y, \quad \lambda \neq \mu >0.
\end{eqnarray*}

\noindent
He computes $M$ and plots the contour lines of $M$ for two cases: $\lambda =1, \, \mu =2$ (shown in the left-hand panel of his figure 4) and $\lambda =2, \, \mu =1$ (shown in the right-hand panel of his figure 4), both using $\tau =10$. In both cases the contour lines are straight lines and fail to reveal the stable and unstable manifolds. However as we have pointed out in previous examples, this is not the correct way to rule out that $M$ detects stable and unstable manifolds.

\begin{figure*}[htbp!]
  \centering
  \subfigure[Contours of $\partial_{x_{0}} M$,  $\lambda =1, \, \mu =2, \tau=10$]{\includegraphics[width=0.48\linewidth]{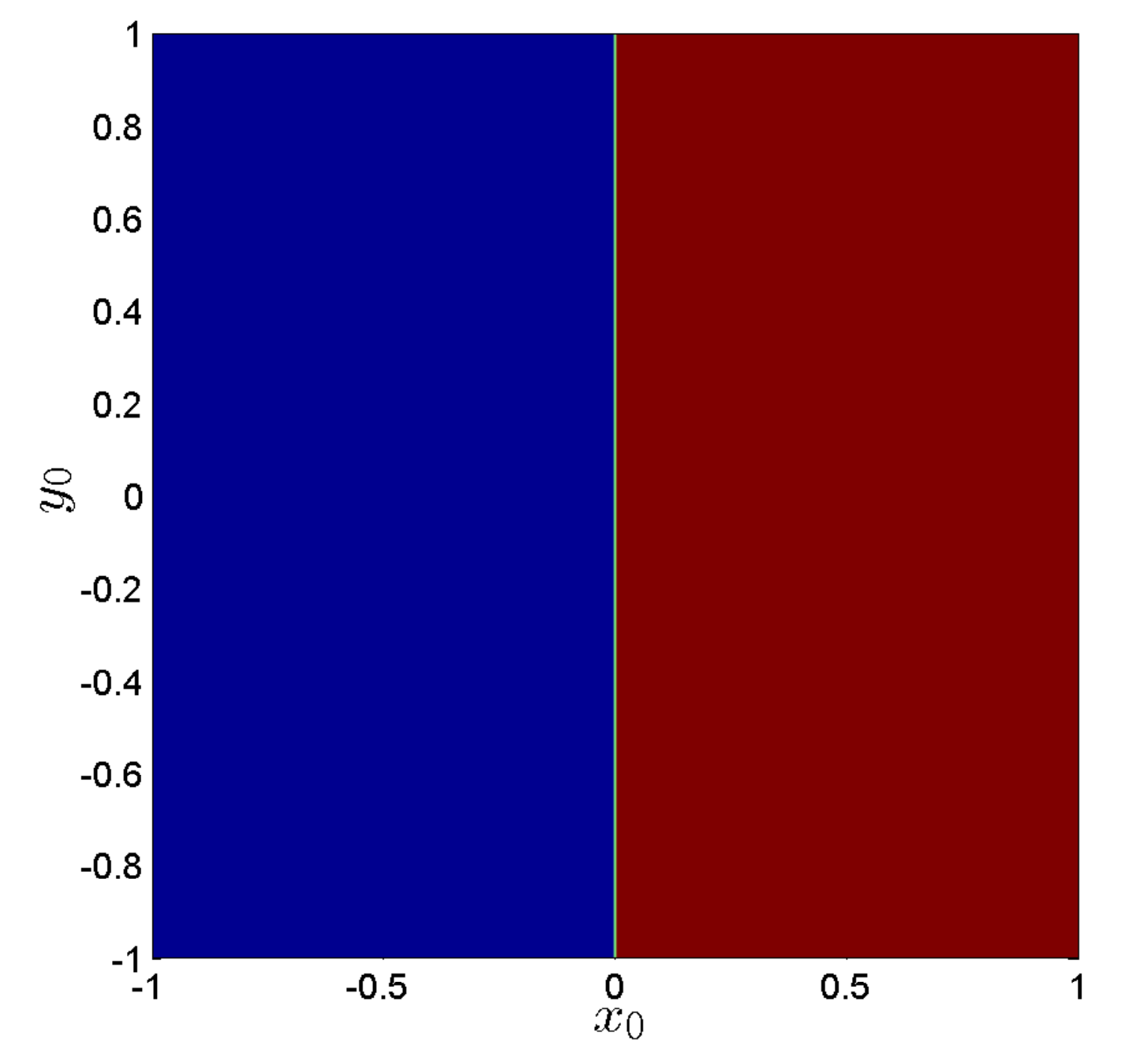}}
  \subfigure[Contours of $\partial_{y_{0}} M$,  $\lambda =1, \, \mu =2, \tau=10$]{\includegraphics[width=0.48\linewidth]{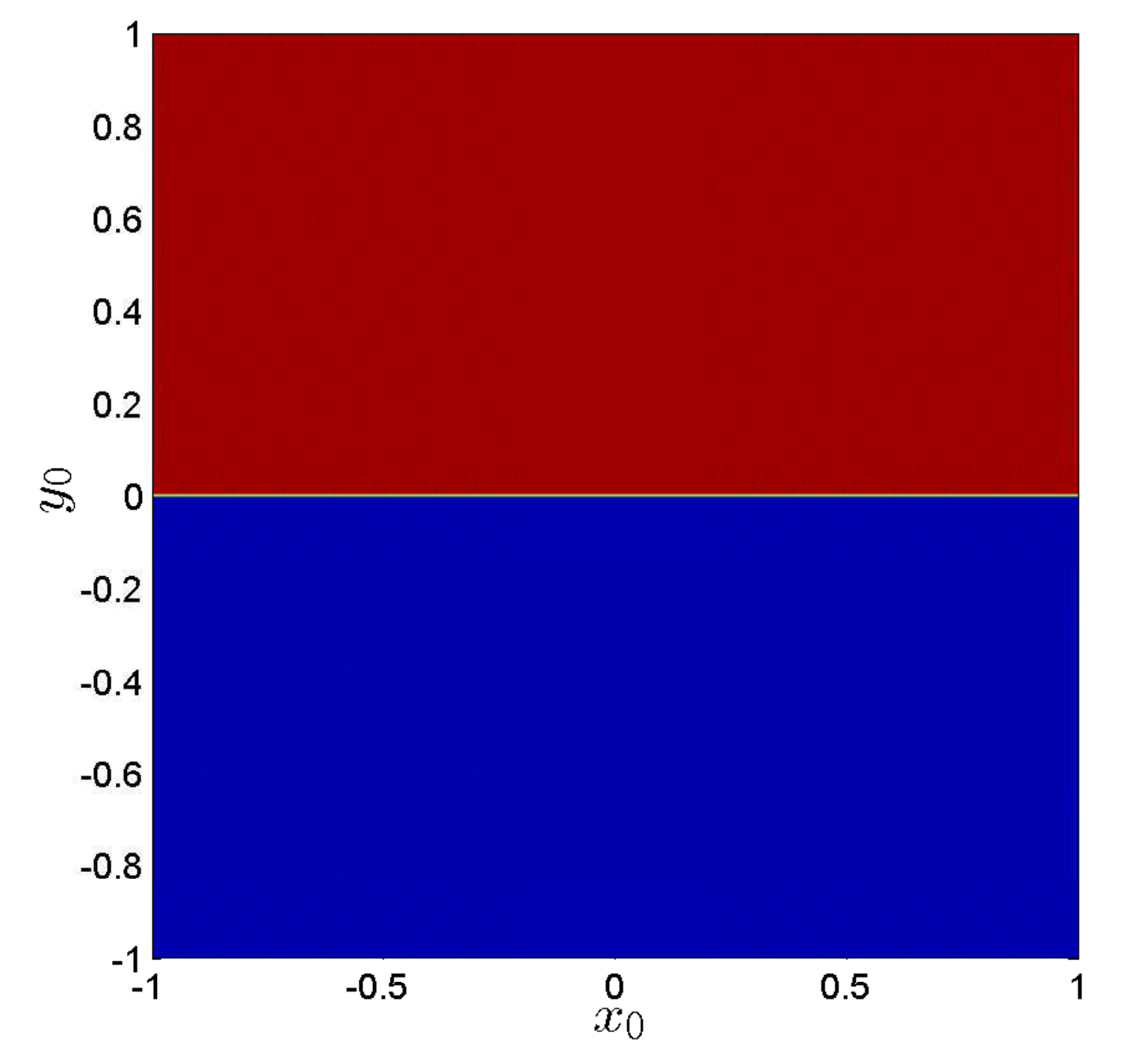}} \\
  \subfigure[Contours of $\partial_{x_{0}} M$,  $\lambda =2, \, \mu =1, \tau=10$]{\includegraphics[width=0.48\linewidth]{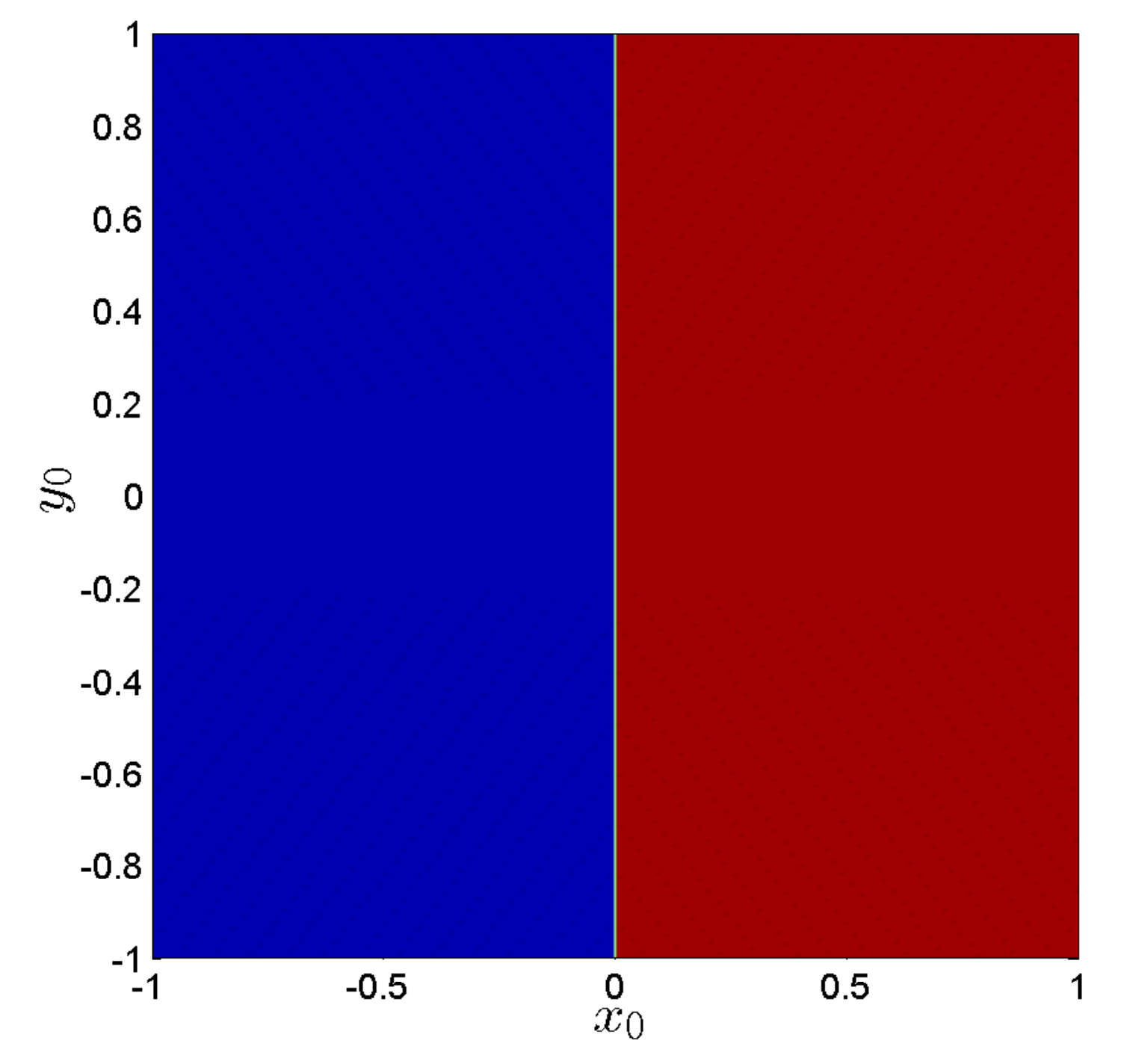}}
  \subfigure[Contours of $\partial_{y_{0}} M$,  $\lambda =2, \, \mu =1, \tau=10$]{\includegraphics[width=0.48\linewidth]{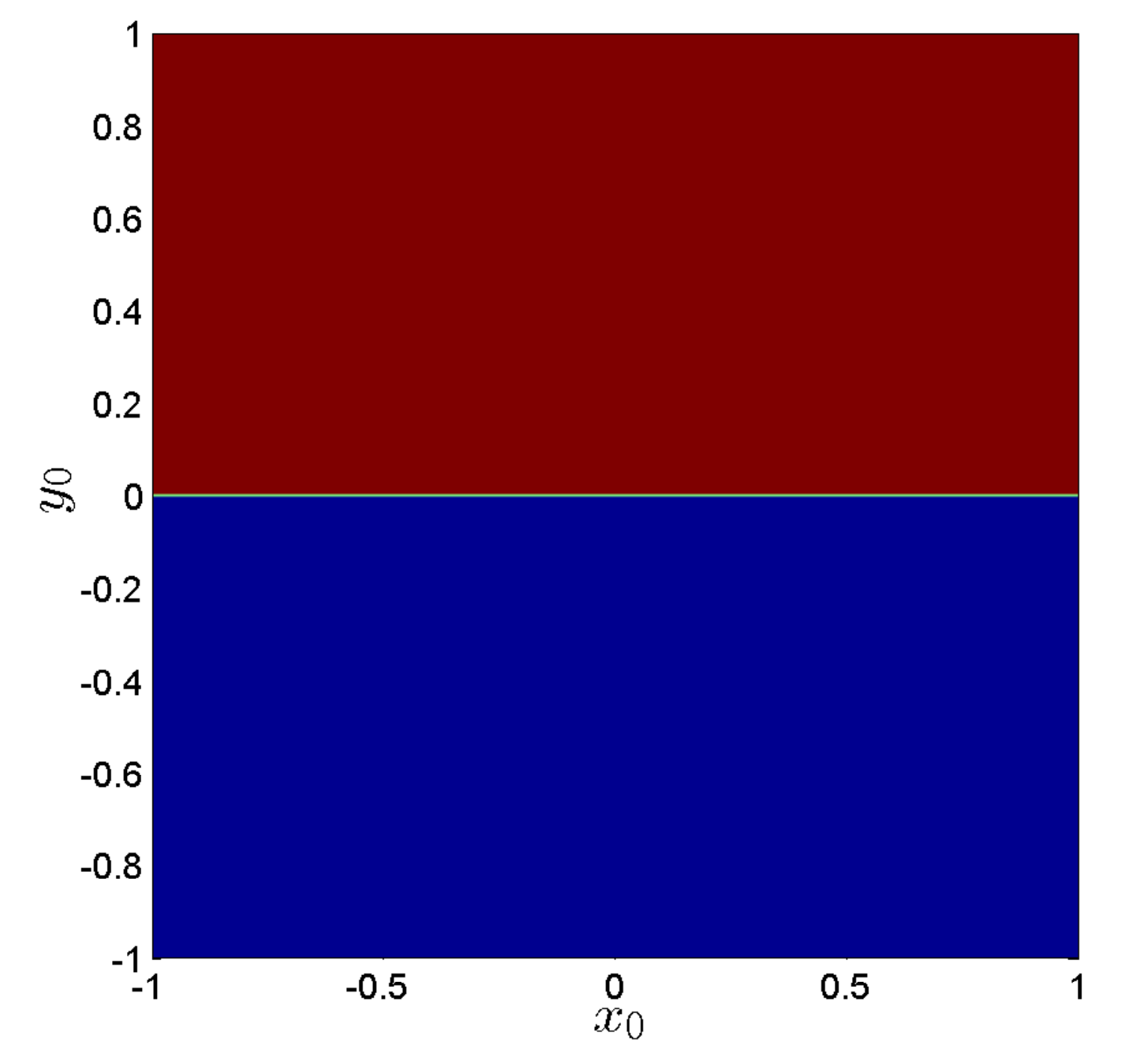}} \\
  \caption{In panels a) and c) the stable manifold is clearly identified and in panels b) and d) the unstable manifold is clearly identified. This figure should be compared with figure 4 of the Comment of Ruiz-Herrera.}
  \label{Afig4}
\end{figure*}

In panels a) and b) we plot contours of the the derivatives $\partial_{x_{0}} M$ and $\partial_{x_{0}} M$, respectively, for the case $\lambda =1, \, \mu =2, \tau=10$. We see that the stable and unstable manifolds are clearly identified in these plots, just as in the earlier examples.
In panels c) and d) we plot contours of the the derivatives $\partial_{x_{0}} M$ and $\partial_{x_{0}} M$, respectively, for the case $\lambda =2, \, \mu =1, \tau=10$. Also here in this case we see that the stable and unstable manifolds are clearly identified in these plots.

Apparently  Ruiz-Herrera's reasoning behind his Example 1 and Example 4 was that differing expansion and contraction rates of the saddle could cause the method of Lagrangian Descriptors to fail. Clearly  this is not true.

\section{Conclusions of Ruiz-Herrera}

Ruiz-Herrera ends his Comment with a section of conclusions based on his ``counterexamples''. We have shown that, due to his lack of understanding of the singular properties of the LDs, all of his ``counterexamples''  and therefore {\em all} of his conclusions  are wrong.  

He ends his conclusions section by commenting on the recent development of LDs for discrete time systems described in \cite{carlos}. He claims that the results of his Comment are applicable to the time 1 map of his examples and consequently this fact implies that the results of \cite{carlos} are incorrect. As we have seen, such an assertion is wrong.

Finally, we note that Ruiz-Herrera has recently published  another collection of ``counterexamples'' to the method of Lagrangian Descriptors (\cite{RH2015}).  These ``counterexamples'' also suffer from an incorrect  understanding of the singular properties of the of LDs, and are also wrong.

\section*{\bf Acknowledgments.} The research of FB-I, JC, VJG-G, CL, AMM and CM is supported by the MINECO under grant MTM2014-56392-R. The research of SW is supported by  ONR Grant No.~N00014-01-1-0769.  JC is grateful to the LABEX Lyon Institute of Origins (ANR-10-LABX-0066) of the Universit\'e de Lyon for its financial support. We acknowledge support from MINECO: ICMAT Severo Ochoa project SEV-2011-0087.

\end{document}